\documentclass[11pt,a4paper]{article}
\usepackage[utf8]{inputenc}
\usepackage{graphicx,verbatim,array,multicol,courier}
\usepackage{amsmath,amssymb,amsthm,mathrsfs, mathtools}
\usepackage{color, graphics, graphicx, siunitx}
\usepackage{amsfonts, dsfont, bm}
\usepackage{natbib} 
\usepackage{geometry}
\usepackage{enumitem,float}
\usepackage{lscape}
\usepackage[pdftex,bookmarks,colorlinks]{hyperref}
\usepackage[dvipsnames]{xcolor}
\usepackage{tabularx,multirow}
\usepackage{siunitx}
\usepackage[labelformat=simple]{subfig}
\usepackage[linesnumbered,ruled,vlined]{algorithm2e}
\usepackage{booktabs}

\usepackage{accents}
\newlength{\dhatheight}

\newlist{inparaenum}{enumerate}{2}
\setlist[inparaenum,1]{label=(\alph*)}
\setlist[inparaenum,2]{label=(\roman{inparaenumi}\emph{\alph*})}

\newcommand\round[1]{\left[#1\right]}

\makeatletter
\def\adl@drawiv#1#2#3{%
        \hskip.5\tabcolsep
        \xleaders#3{#2.5\@tempdimb #1{1}#2.5\@tempdimb}%
                #2\z@ plus1fil minus1fil\relax
        \hskip.5\tabcolsep}
\newcommand{\cdashlinelr}[1]{%
  \noalign{\vskip\aboverulesep
           \global\let\@dashdrawstore\adl@draw
           \global\let\adl@draw\adl@drawiv}
  \cdashline{#1}
  \noalign{\global\let\adl@draw\@dashdrawstore
           \vskip\belowrulesep}}
\makeatother

\numberwithin{equation}{section}
\theoremstyle{definition}
\newtheorem{defi}{Definition}[section]
\newtheorem{cond}[defi]{Condition}
\theoremstyle{plain}
\newtheorem{theo}[defi]{Theorem}
\newtheorem{prop}[defi]{Proposition}

\theoremstyle{remark}

\theoremstyle{example}

\newcommand{\diff}{\mathrm{d}}

\definecolor{navy}{rgb}{0,0,0.502}
\definecolor{brown}{rgb}{0.59, 0.29, 0.0}
\def\indic{\mathds{1}}

\hypersetup{colorlinks,%
	citecolor=blue,%
	filecolor=green,%
	linkcolor=red,%
	urlcolor=violet,%
}

\newcommand{\Real}{\mathbb{R}}

\newcommand{\Expect}{\mathbb{E}}

\newcommand{\RV}{\mathcal{RV}}
\newcommand{\Prob}{\mathbb{P}}
\newcommand{\Eal}{\mathscr{E}}

\newcommand{\bfX}{{\boldsymbol{X}}}
\newcommand{\bfW}{{\boldsymbol{W}}}
\newcommand{\bfU}{{\boldsymbol{U}}}

\newcommand{\bfzero}{{\boldsymbol{0}}}
\newcommand{\bfone}{{\boldsymbol{1}}}
\newcommand{\bfx}{{\boldsymbol{x}}}
\newcommand{\bfw}{{\boldsymbol{w}}}
\newcommand{\bfu}{{\boldsymbol{u}}}
\newcommand{\simp}{{\mathcal{S}}}

\title{Marginal expected shortfall inference under multivariate regular variation}
\author{S. A. Padoan,  S. Rizzelli and M. Schiavone}

\begin{document}

\maketitle
\begin{abstract}
Marginal expected shortfall is unquestionably one of the most popular systemic risk measures. Studying its extreme behaviour is particularly relevant for risk protection against severe global financial market downturns. In this context, results of statistical inference rely on the bivariate extreme values approach, disregarding the extremal dependence among a large number of financial institutions that make up the market. In order to take it into account we propose an inferential procedure based on the multivariate regular variation theory. We derive an approximating formula for the extreme marginal expected shortfall and obtain from it an estimator and its bias-corrected version. Then, we show their asymptotic normality, which allows in turn the confidence intervals derivation. Simulations show that the new estimators  greatly improve upon the performance of existing ones and confidence intervals are very accurate. An application to financial returns shows the utility of the proposed inferential procedure. Statistical results are extended to a general $\beta$-mixing
context that allows to work with popular time series models with heavy-tailed innovations.
\end{abstract}

%
\section{Introduction}\label{sec:intro}
%

Systemic risk measures are of crucial importance in assessing the risk of the global financial system in today's connected world. A simple example of a coherent measure is the Expected Shortfall (ES) of the system, i.e. the expected loss based on the financial market return, given that a systemic shortfall event has been experienced \citep[e.g.][]{acharya2017measuring}.
%
%
From a statistical viewpoint, let $X_j$ be a random variable representing the negative return (referred in the sequel as returns for brevity) for the $j$th financial institution (bank), with distribution $F_j$ for $j\in\{1,\ldots,d\}$, and $R=X_1+\ldots+X_d$ be the return of the market. Then, the ES of the financial system is 
$$
\Theta(\tau)=\Expect(R| R> Q_R(\tau)),
$$
where for a random variable $X$, $Q_X(\tau)=U_X(1/(1-\tau))$ is the $\tau$-quantile of $F_X$ also known as Value-at-Risk with confidence level $\tau$, with  $\tau=1-p$ and $p\in(0,1)$ is meant as a small exceeding probability and where $U_X(t)=\inf\{x\in\Real: F_X\geq 1-t^{-1}\}$, for $t>1$. The Marginal Expected Shortfall (MES) is a quantity closely related to the ES of the financial market and it plays a pivotal role in the systemic risk measurement procedure. The MES is the expected equity loss of an individual bank given indeed the occurrence of the market downturn
%
\citep[e.g.][]{acharya2017measuring}. It is defined as
\begin{equation}\label{eq:MES}
\theta_{j}(\tau)=\Expect\left(X_j\bigg{\vert} R> Q_R(\tau)\right), \quad j\in\{1,\ldots,d\},
\end{equation}
and it satisfies the property that $\Theta(\tau)=\sum_{j\leq d}\theta_j(\tau)$, for any $\tau\in(0,1)$. This makes it a particularly useful tool to be employed by regulatory policies to establish the contribution of individual institutions to the risk of the entire financial system, thus requesting, for example, their capital increase \citep[e.g.][]{tarashev2010}.

In real financial analyses, assessing the extreme MES (corresponding to confidence level $\tau$ approaching one) is of crucial importance in oder to understand the risk of severe financial market downturns. Its proper computation requires the quantification of the dependence level among the components of $\bfX$. In this regard, we study the MES behaviour under the multivariate regular variation framework  \citep[e.g.][]{resnick2007}.

Our first contribution is to provide a simple closed form expression for approximating formula \eqref{eq:MES}, when $\tau\uparrow1$. This result agrees with the finding in the paper \citet{joe2011tail}, whose authors work in more generality but without providing an explicit formula that can be easily used for the statistical estimation of the MES. We determine the speed with which our proposal approximates the MES, when $\tau\uparrow 1$, under some standard multivariate regular variation conditions at the density level. Several estimators of the extreme MES have been proposed in the literature, e.g., \citet{cai2015}, \citet{cai2017}, \citet{di2018}, \citet{goegebeur2021}, to name a few; see also \citet{daouia2018} for an expectile based version. However, all of them rely on the bivariate extreme value theory rather than a proper multivariate theory in arbitrary dimension. Our second contribution is to provide a MES estimator grounded on our approximate formula and  we establish its asymptotic normality. Through a simulation study, we show the superior performance of our estimator compared to the two estimators that are proposed in \citet{cai2015}, which are a benchmark in the standard literature. The gain of our proposal over the competitors is particularly relevant in high dimensions. We propose a further bias-corrected version of it and we establish again the asymptotic normality. We show that it outperforms the competitors and behaves steadily for a large interval of the effective sample fraction, between $2\%$ up to $30\%$, and this result is particularly useful as it simplifies the notorious problem of selecting the latter value in real applications. Underestimation of the variance term in the limiting normal distribution heavily affects the accuracy of the confidence interval that is constructed on its basis. We propose then an adjusted version of the confidence intervals initially derived, whose coverage probability is close to the nominal level for a large range of the effective sample fraction. Most of the existing results on the behaviour of the extreme MES are derived to work with independent data, apart from the recent work \citet{davison2022tail}.  As a final contribution, we extend the asymptotic results of our estimator to the case where it is constructed from a $\beta$-mixing heavy-tailed strictly stationary multivariate time series. Our methods are integrated of the {\tt R} package {\tt ExtremeRisks} available on CRAN.

The paper is organised as follows. Section \ref{sec:background} provides the statical background and  our main approximating formula for the MES. Section \ref{sec:newapproach} discusses the MES estimation and describe a finite-sample procedure for constructing confidence intervals. The methods are tested by a simulation study in Section \ref{sec:simulation} and applied to financial data in Section \ref{sec:realdata}. The paper ends discussing the extension of the theoretical results in a time series context in Section \ref{sec:extension}. The online supplementary material contains proofs and further results on the real data analysis.

%
\section{Statistical model}\label{sec:background}
Let $\bfX=(X_1,\ldots,X_d)^\top$ be a $d$-dimensional random vector on $[0,\infty)^d$.  Consider the transformation to radial $R=\|\bfX\|$ and 
angular $\bfW=R^{-1}\bfX$ components,  where $\|\cdot\|$ is any norm in $[0,\infty)^d$. Let $\simp=\{\bfw\in[0,\infty)^d:\|w\|=1\}$ be the unit sphere in $[0,\infty)^d$ with respect to $\|\cdot\|$. We assume that the distribution $F$ of $\bfX$ is regularly varying with index $1/\gamma$ and $\gamma>0$, in symbols $F\in\RV(1/\gamma)$, i.e. there is a Radon measure $\nu$ on $\Eal:=[\boldsymbol{0}, \boldsymbol{\infty}]\setminus\{\boldsymbol{0}\}$  such that for any relative compact Borel set $B\subset\Eal$ 
with $\nu(\partial B)=0$ we have 
\begin{equation}\label{eq:exponent_measure*}
\frac{\Prob(t^{-1}\bfX\in B)}{\Prob(R>t)} \to \nu(B), \quad t\to\infty
\end{equation}
%
and it satisfies the homogeneity property  $\nu(xB)=x^{-1/\gamma}\nu(B)$ for any Borel set $B\subset \Eal$ and any $x>0$ 
\citep[e.g.][Ch. 6] {jessen2006regularly, resnick2007}. Result in \eqref{eq:exponent_measure*} is equivalent to having that $\Prob(R>x)$ is regularly varying with index $1/\gamma$ and
$$
\Prob(\bfW\in B| R>t)\to H(B), \quad t\to\infty,
$$
for any Borel set $B\subset\simp$ 
with $H(\partial B)=0$, where 
$H$ is a probability measure on $\simp$. We denote the marginal angular average as
$$
\overline{w}_j=\int_\simp w_jH(\diff \bfw).
$$
The above two measures are related as follows, 
$$
\lim_{t\to\infty}\frac{1-F(t\bfx)}{1-F(t\bfone)}=\frac{\nu\left([\bfzero,\bfx]^\complement\right)}{\nu\left([0,\bfone]^\complement\right)}=\frac{\gamma^{-1}}{\nu\left([0,\bfone]^\complement\right)}\int_{\simp}\max_{1\leq j\leq d}\left(\frac{w_j}{x_j}\right)^{1/\gamma}H(\diff\bfw)
$$
for all $\bfx\in\Eal$. In the case that $\gamma=1$, then $\nu$ and $H$ are the canonical exponent and angular probability measures, respectively, e.g.  in \citet[][Theorem 6.1]{resnick2007}, \citet[][Ch. 6.1.3]{dehaan+f06} and \citet[][Ch. 4.2]{falk2010}.

We are now ready to state our first main probabilistic result. In the sequel we assume for simplicity that $\|\cdot\|$ in the above framework is the $ L_1$-norm, denoted as $\|\cdot\|_1$.
\begin{prop}\label{prop:limprob}
Let $\bfX$ be a random vector whose distribution satisfies $F\in\RV(1/\gamma)$, $\gamma>0$. For $1\leq i,j\leq d$ assume also
$$
\lim_{x \to \infty}\frac{\Prob(X_i>x)}{\Prob(X_j>x)}=c_{i,j}\in (0,\infty).
$$
Then, for any $j\in\{1,\ldots,d\}$ we have 
$$
\lim_{\tau\to 1^{-}}\frac{\theta_j(\tau)}{Q_R(\tau)} = \frac{\overline{w}_j}{1-\gamma}.
$$
\end{prop}
By this result we obtain the following approximation
$$
\theta_j(\tau)\approx Q_R(\tau)\frac{\overline{w}_j}{1-\gamma},\quad \tau\uparrow 1.
$$
Next theorem quantifies how fast the above convergence result is,  further assuming second-order conditions at the density level. 
In particular, assume that $F$ allows for a density $f$ and $H$ allows for a bounded, continuous density $h:\simp\to(0,\infty)$ such that
%
%
\begin{equation*}\label{eq:uniform_limit}
\lim_{t\to\infty} \sup_{\bfw\in\simp}\left|\frac{t^{d}f(t\bfw)}{\Prob(R>t)}-\frac{h(\bfw)}{\gamma}\right|=0.
\end{equation*}
As a result, $\nu$ also allows for a density function $q$ that satisfies $q(t\bfx)=t^{-(d+1/\gamma)}q(\bfx)$, for all $t>0$ with $\bfx\in\Eal$, and for all sets $B$ in the Borel $\sigma$-field of $\Eal$ we have
$$
\nu(B) = \int_{\bfx\in B}q(\bfx)\diff \bfx = \frac{1}{\gamma}\int_{\{r>0, \bfw\in\simp:\bfx\in B\}} r^{-1/\gamma-1} h(\bfw) \diff r \diff \bfw.
$$
%
Next are the additional  second-order conditions required.
\begin{cond}\label{cond:mvt_sec_order}
For some $\rho<0$ and any $j\in\{1,\ldots,d\}$: 
\begin{inparaenum}
\item \label{cond. Hall_class} there is a rate function $A_{X_j}(t)$ such that $\lim_{t\to\infty} t^{-\rho}A_{X_j}(t) \in \Real \setminus \{0\}$ and 
$$
\lim_{t\to\infty}\frac{\frac{1-F_{X_j}(tx)}{1-F_{X_j}(t)}-x^{-1/\gamma}}{A_{X_j}(t)}=x^{-1/\gamma} \frac{x^\rho-1}{\rho};
$$
\item \label{cond: mvt_dep} there is a nonnull function $\lambda_j$, integrable on sets bounded away from zero,  such that  $\sup_{\bfw\in\simp}\vert \lambda_j(\bfw)\vert<\infty$, and a function $q_j(\bfx)\propto r^{-d-1/\gamma} h(\bfw)$ such that 
$$
\lim_{t\to\infty}
\frac{\frac{t^df(t\bfx)}{\Prob(X_j>t)}-q_j(\bfx)}{A_{X_j}(t)}=\lambda_j(\bfx),\,\forall\, \bfx\in(0,\infty)^d,
$$
and
$$
\lim_{t\to\infty}\sup_{\bfw\in\simp}\left\vert
\frac{\frac{t^df(t\bfw)}{\Prob(X_j>t)}-q_j(\bfw)}{A_{X_j}(t)}
-\lambda_j(\bfw)
\right\vert=0.
$$
\end{inparaenum}
\end{cond}

\begin{theo}\label{theo:limprob_speed}
If Condition \ref{cond:mvt_sec_order} is satisfied, then for each $j\in\{1,\ldots,d\}$ there is a constant $c_j\in\Real$ such that for $\tau=1-1/t$,
\begin{equation*}
A_{X_j}^{-1}(U_{X_j}(t))\left(\frac{\theta_j(\tau)}{Q_R(\tau)\frac{\overline{w}_j}{1-\gamma}}-1\right)\to c_j, \quad t\to\infty.
%
\end{equation*}
\end{theo}
%


%
\section{Inference}\label{sec:newapproach}

Let $\bfX_1,\ldots,\bfX_n$ be a sample of i.i.d. copies of $\bfX$ whose distribution satisfies $F\in\RV(1/\gamma)$ with $\gamma>0$. The aim is to estimate $\theta_j(\tau)$, for any $j\in\{1,\ldots,d\}$, corresponding to an extreme level $\tau_n$  such that $\tau_n\to 1$ as $n\to\infty$ and in particular $n(1-\tau_n)\to a>0$ as $n\to\infty$. 
This task is achieved combining estimates of $\gamma$ and $Q_R$ suitably, which are obtained from the so-called effective sample fraction $k_n/n$ of the data. In particular, we assume $k_n\to\infty$ as $n\to\infty$ and $k_n=o(n)$.

Our proposal to estimate  $\theta_j(\tau_n)$ is as follows. First, we estimate the tail index $\gamma$ by the Hill estimator \citep{hill1975},
$$
\widehat{\gamma}_n=\sum_{i=1}^{k_n} \frac{i}{k_n} \ln \frac{R_{(n-i+1,n)}}{R_{(n-i,n)}},
$$
where $R_{(1,n)}\leq\cdots\leq R_{(n,n)}$ are the $n$-order statistics of the sample $R_1,\dots,R_n$. Let 
$$
\widehat{Q}_R(\tau_n)=R_{(n-k_n,n)}\left(\frac{n(1-\tau_n)}{k_n}\right)^{-\widehat{\gamma}_n},
$$
be the Weissman estimator \citep{weissman1978} of radial quantile at the level $\tau_n$ and
$$
\widehat{\overline{w}}_{j,n}=\frac{1}{k_n} \sum_{i=1}^{n}\frac{X_{i,j}}{R_i}
\indic \left( R_i > R_{(n-k_n,n)}
\right)
$$
be the estimator of the mean of the $j$th angular component at the intermediate level $1-k_n/n$. Finally, an estimator of $\theta_j(\tau_n)$ is 
$$
\widehat{\theta}_j(\tau_n)=\widehat{Q}_R(\tau_n)\frac{\widehat{\overline{w}}_{j,n}}{1-\widehat{\gamma}_n}.
$$

Under Condition \ref{cond:mvt_sec_order}, the distribution of the radius $F_R$ belongs to the so-called Hall-Welsh class \citep{hall1985adaptive}, 
i.e. it allows the representation
$$
U_R(t)=C_R t^{\gamma}\left(
1+\frac{A_R(t)}{\rho_R}+o(t^{\rho_R})\right),\quad A_R(t)=\beta_R t^{\rho_R}, 
$$
where $A_R$ and $\rho_R$ are the rate function and second-order parameter of the radial component, $C_R>0$, $\beta_R \neq 0$ and $\rho_R=\gamma \rho$, see Lemma B.6 in the supplementary material. Next, we consider the following adjusted or biased-corrected version of the tail index estimator and Weissman estimator of radial quantile,
\begin{align*}
	&\widehat{\gamma}_n^{\text{(Adj)}} := \widehat{\gamma}_n\left(
	1-\frac{\widehat{\beta}_{R,n}}{1-\widehat{\rho}_{R,n}}\left(\frac{n}{k_n}\right)^{\widehat{\rho}_{R,n}}
	\right),\\
	&\widehat{Q}_R^{\text{(Adj)}}(\tau):= R_{(n-k_n,n)}
	\left(\frac{n(1-\tau_n)}{k_n}
	\right)^{-\widehat{\gamma}_n^{\text{(Adj)}}}\exp\left(C_\tau(n;\widehat{\beta}_{R,n},\	\widehat{\rho}_{R,n})\right),
\end{align*}
where $\widehat{\beta}_{R,n}\equiv \widehat{\beta}_R(s_n)$ and $\widehat{\rho}_{R,n}\equiv \widehat{\rho}_R(s_n)$ are the estimators of $\beta_R$ and $\rho_R$ in Sections 2.1 and 2.2 of \cite{gomes2007}, with $s_n$ that is another intermediate sequence different from $k_n$ such that $s_n\to\infty$ as $n\to\infty$ and $s_n=o(n)$ 
and
\begin{align*}
	C_\tau(n;\widehat{\beta}_{R,n},\widehat{\rho}_{R,n}):= \widehat{\beta}_{R,n}\left(
	\frac{n}{k_n}
	\right)^{\widehat{\rho}_{R,n}}
	\frac{(k_n/(n(1-\tau_n)))^{\widehat{\rho}_{R,n}}-1}{\widehat{\rho}_{R,n}}.
\end{align*}
On the basis of these refined estimators we can now define the following adjusted estimator of $\theta_j(\tau)$, 
$$
\widehat{\theta}_{j}^{\text{(Adj)}}(\tau_n):=\widehat{Q}_{R,n}^{\text{(Adj)}}
\left({\tau_n}\right)\frac{\widehat{\overline{w}}_{j,n}}{1-\widehat{\gamma }_n^{\text{(Adj)}}},
$$
for any $j\in\{1,\ldots,d\}$. The unadjusted and adjusted versions of the MES estimator satisfy both the asymptotic normality property.
%
%
\begin{theo}\label{theo:est}
Assume that Condition \ref{cond:mvt_sec_order} is satisfied. Assume also that $k_n=o(n)$, 
$\sqrt{k_n}A_R(n/k_n)\to\lambda_R\in \Real$ as $n\to\infty$. 
Then, for any $j\in\{1,\ldots,d\}$, as $n\to\infty$
%
$$
\frac{\sqrt{k_n}}{\log(k_n/n(1-\tau_n))}\log \frac{\widehat{\theta}_j(\tau_n)}{\theta_j(\tau_n)}
	\stackrel{d}{\to}
	\mathcal{N}\left(\frac{\lambda_R}{1- \rho_R},\gamma^2\right).
$$
Furthermore, if $s_n=o(n)$ and $\sqrt{s_n}A_R(n/s_n)\to \infty$ as $n\to\infty$, then as $n\to\infty$
$$
\frac{\sqrt{k_n}}{\log(k_n/n(1-\tau_n))}\log \frac{\widehat{\theta}^{\text{(Adj)}}_j(\tau_n)}{\theta_j(\tau_n)}
	\stackrel{d}{\to}
	\mathcal{N}\left(0,\gamma^2\right).
$$
\end{theo}
Exploiting the first result we have that an approximate $(1-\alpha)100\%$ confidence interval for $\theta_j(\tau_n)$ can be obtained as
\begin{equation}\label{eq:app_ci_mse}
\left[
\widehat{\theta}_j(\tau_n)\left(\frac{n(1-\tau_n)}{k_n}\right)^{\widehat{b}_n+ z_{1-\alpha/2}\widehat{\gamma}_n/\sqrt{k_n}};\;
\widehat{\theta}_j(\tau_n)\left(\frac{n(1-\tau_n)}{k_n}\right)^{\widehat{b}_n-z_{1-\alpha/2}\widehat{\gamma}_n \sqrt{k_n}}
\right],
\end{equation}
where the bias term $\lambda_R/(1-\rho_R)$ is estimated by $\sqrt{k_n}\widehat{b}_n$, where $\widehat{b}_n=\widehat{\gamma}_n \widehat{\beta}_{R,n}(n/k_n)^{\widehat{\rho}_{R,n}}/(1-\widehat{\rho}_{R,n})$ and where $z_p$ is the $p$-quantile of the standard normal distribution. In practice, the observed values of $ \widehat{\beta}_{R,n}$ and $\widehat{\rho}_{R,n}$ are computed using the {\tt R} package {\tt evt0} \citep{manjunath2013evt0}. Alternatively, one can consider the second asymptotic result and end-up with a similar interval but where $\widehat{\theta}_j(\tau_n)$ is replaced by $\widehat{\theta}^{\text{(Adj)}}_j(\tau_n)$ and $\widehat{b}_n$ by zero.  In practice, both intervals perform similarly, although in particular the adjusted version does better for large values of $k_n$. Underestimation of the variance term $\gamma^2$ in Theorem \ref{theo:est} entails a poor coverage probability of the aforementioned intervals. In this regard to improve the interval's performance one can consider a further adjustment as for example that along the lines of the suggestion on page 287 of \cite{gomes2007}. Since such an approach does not provide the desired results in this context, we propose  alternative bias and variance corrected intervals.  Our new adjustment is a slightly more refined version of the confidence interval \eqref{eq:app_ci_mse}, obtained on the basis of the following simple considerations. First, we have verified that the element $\widehat{\overline{w}}_{j,n}$ contributes little to the total variability of the estimator $\widehat{\theta}_{j}(\tau_n)$ in contrast to the other two components which are instead much more relevant contributors. We then focus on the following expansion
\begin{eqnarray*}
\log \frac{\widehat{\theta}_j(\tau_n)}{\theta_j(\tau_n)}\frac{\overline{w}_j}{\widehat{\overline{w}}_j}&=& (\widehat{\gamma}_n-\gamma)\left(\frac{\sqrt{k_n}}{c_n}+\frac{1}{1-\gamma}\right)+\log\frac{R_{n-k_n+1,n}}{U(n/k_n)}  + \text{O}_p\left((\gamma-\widehat{\gamma}_n)^2\right)\\
&=&  (\widehat{\gamma}_n-\gamma)\left(\log\frac{\sqrt{k_n}}{c_n}+\frac{1}{1-\gamma}\right)+\frac{R_{n-k_n+1,n}}{U(n/k_n)}-1 \\
&+&\text{O}_p\left(\left(\frac{R_{n-k_n+1,n}}{U(n/k_n)}-1\right)^2\right)  + \text{O}_p\left((\gamma-\widehat{\gamma}_n)^2\right)\\
&=& (\widehat{\gamma}_n-\gamma)\left(\frac{\sqrt{k_n}}{c_n}+\frac{1}{1-\gamma}\right)+\frac{R_{n-k_n+1,n}}{U(n/k_n)}-1 + \text{O}_p\left(1/k_n\right).
\end{eqnarray*}
where $c_n=\sqrt{k_n}/(\log(k/n(1-\tau_n)))$. On this basis we obtain that (approximately)
\begin{eqnarray*}
\Expect\left(c_n\log\frac{\widehat{\theta}_j(\tau_n)}{\theta_j(\tau_n)}\frac{\overline{w}_j}{\widehat{\overline{w}}_j} \right)&=&\frac{\lambda_R}{1-\rho_R}\left(1+\frac{c_n}{\sqrt{k_n}}\frac{1}{1-\gamma}\right),\\
\text{Var}\left(c_n\log\frac{\widehat{\theta}_j(\tau_n)}{\theta_j(\tau_n)}\frac{\overline{w}_j}{\widehat{\overline{w}}_j} \right)&=&\gamma^2\left(1+\frac{2}{1-\gamma}\frac{c_n}{\sqrt{k_n}}+2\frac{c_n^2}{k_n}\right)
\end{eqnarray*}
and we then propose the refined approximate $(1-\alpha)100\%$ confidence interval for $\theta_j(\tau_n)$,  
\begin{equation}\label{eq:app_ci_corrected_mse}
\left[
\widehat{\theta}_j(\tau_n)\left(\frac{n(1-\tau_n)}{k_n}\right)^{\widehat{b}^*_n+ z_{1-\alpha/2}\widehat{v}_n/\sqrt{k_n}};\;
\widehat{\theta}_j(\tau_n)\left(\frac{n(1-\tau_n)}{k_n}\right)^{\widehat{b}^*_n-z_{1-\alpha/2}\widehat{v}_n \sqrt{k_n}}
\right],
\end{equation}
where 
$$
\widehat{b}^*_n=\widehat{b}_n\left(1+\frac{c_n}{\sqrt{k_n}}\frac{1}{1-\widehat{\gamma}_n}\right),\quad
\widehat{v}_n=\widehat{\gamma}_n\left(1+\frac{2}{1-\widehat{\gamma}_n}\frac{c_n}{\sqrt{k_n}}+2\frac{c_n^2}{k_n}\right)^{1/2}.
$$
An alternative to this last interval is obtained by substituting $\widehat{\theta}_j(\tau_n)$ by $\widehat{\theta}^{\text{(Adj)}}_j(\tau_n)$ and $\widehat{b}^*_n$ by zero.

%

%
\section{Simulation study}\label{sec:simulation}
%
Through an extensive simulation study, we show the finite samples performance of our proposed estimators and those of their corresponding confidence intervals. The simulation experiment is based on: three different bi-variate distributions, one four-variate distribution and one fifteen-variate distribution. We consider 3 Archimedean copulae families and an elliptical one. We recall that an Archimedean copula with generator $\varphi$ is a distribution defined on the $d$-dimensional hypercube as
$$
C(\bfu)=\varphi^{-1}(\varphi(u_1)+\cdots+\varphi(u_d)),\quad \bfu\in[0,1]^d,
$$
where $\varphi:(0,1]\to[0,\infty)$ is a convex and strictly decreasing
function with $\varphi(1)=0$ and $\varphi(t)\to\infty$ as $t\to 0$ \citep[see][for details]{joe2014}.
We focus on the following classes. The Clayton, with generator $\varphi(u)=\delta^{-1}(u^{-\delta}-1)$ for $\delta>0$, with $\delta\to0$ and $\delta\to\infty$ representing the case of independence and complete dependence among the $\bfu$ components. The Gumbel, with generator $\varphi(u)=(-\log(u))^\delta$ for $\delta\geq 1$, with $\delta=1$ and $\delta\to\infty$ representing the case of independence and complete dependence. The Joe \citep[][Ch. 4.7]{joe2014}, with generator $\varphi(u)=-\log(1-(1-u)^{\delta})$ for $0\leq \delta\leq 1$, with $\delta=1$ and $\delta\to 0$ representing the case of independence and complete dependence.

In the simulations we consider the following specific models.
\begin{enumerate}
\item[(i)]  [Clayton-Half-$t$ model] Let $\bfU$ follow a Clayton copula with dependence parameter $\delta=3$. Take $X_j=F^{\leftarrow}((U_j+1)/2)$ with $j=1,2$, where $F$ is a Student-$t$ distribution with $5/2$ degrees of freedom. Then, the marginal distributions of $\bfX$ are half-$t$ with tail index $\gamma=2/5$.
\item[(ii)] [Gumbel-Burr model] Let $\bfU$ follow a Gumbel copula with dependence parameter $\delta=0.8$. Take $X_j=F^{\leftarrow}(U_j)$ with $j=1,2$, where $F$ is a Burr distribution with parameters equal to $\sqrt{3}$ and $\sqrt{3}$. Then, the marginal distributions of $\bfX$ are Burr with tail index $\gamma=1/3$.
%
%
\item[(iii)] [Student-$t$-Burr model] Let $\bfU$ follow a Student-$t$ copula with dependence parameter $\omega=0.8$. Take $X_j=F^{\leftarrow}(U_j)$ with $j=1,2$, where $F$ is a Burr distribution with parameters equal to $2$ and $2$. Then, the marginal distributions of $\bfX$ are Burr with tail index $\gamma=1/4$.
\item[(iv)] [Gumbel-Half-$t$-Burr-Fr\'echet-Pareto model] Let $\bfU$ follow a Gumbel copula with dependence parameter $\delta=0.7$. Take $X_j=F_j^{\leftarrow}(U_j)$ with $j=1,\ldots,4$, where $F_j$ with $j=1,\ldots,4$ is equal to Half-$t$, Burr, Fr\'echet and Pareto, respectively, all of them with tail index $\gamma=1/5$.
\item[(v)] [Student-$t$-Half-$t$] Let $\bfU$ follow a Student-$t$ copula with $d=15$ and $4$ degrees of freedom and dependence parameters equal to $1$ along the diagonal and equal to $\omega_{i,j}=0.4$ off-diagonal with $i=2,\ldots,15$ and $j>i$. Take $X_j=F_j^{\leftarrow}((U_j+1)/2)$ with $j=1,\ldots,15$, where $F_j$ is a Student-$t$ distribution with $4$ degree of freedom. Then, the marginal distributions of $\bfX$ are half-$t$ with tail index $\gamma=1/4$.
\end{enumerate}
The components of $\bfX$ are asymptotically independent with model (i) and asymptotically dependent for all the other models, i.e. $\Prob(F_j(X_j)>u|F_i(X_i)>u)\to 0$ as $u\to 1$ for all $1\leq i \neq j \leq d$ in the former case, while such a limit is positive in the latter cases. From each model we simulate $M=50,000$ samples of sample size $n=500$  and compute an estimate of $\theta_j(\tau_n)$, where $\tau_n=1-1/500=0.998$, using our estimators $\widehat{\theta}_{j}(0.998)$ and $\widehat{\theta}_{j}^{\text{Adj}}(0.998)$. In particular, with models (i)-(iii) and (v) we compute an estimate of the MES for the component corresponding to $j=1$ and those relative $j=1,\ldots,4$ with the model (iv). The true values of $\theta_{j}(0.998)$ 
are unavailable in closed-form and we have therefore computed them by Monte Carlo simulations (see Table \ref{tab:true_val_MES}). For comparison purposes we also estimate the MES with the estimators
\begin{eqnarray*}
\widehat{\theta}^{\text{(Emp)}} (\tau_n)&=&\left(\frac{k}{n(1-\tau_n)}\right)^{\widehat{\gamma}_n}\frac{1}{k}\sum_{i=1}^n X_1 \mathbb{I}(R_i>R_{(n-k),n}),\\
\widehat{\theta}^{\text{(Cai)}}(\tau_n) &=& \left(\frac{k}{n(1-\tau_n)}\right)^{\widehat{\gamma}_n}X_{(n-k)}\frac{1}{k}\sum_{i=1}^n\mathbb{I}(R_i>R_{(n-k),n}) \left\{\frac{n-\text{rank}(X_i)+1}{k}\right\}^{-\widehat{\gamma}_n}, \label{Cai_Est}
\end{eqnarray*}
proposed in \cite{cai2015} and we refer to them as the competitors.
With the obtained estimates we compute a Monte Carlo approximation of the squared-bias, variance and Mean Squared Error (MSE) of all the estimators. We also compute a Monte Carlo approximation for the coverage probability of the confidence intervals discussed in the previous section, with $95\%$ nominal level.
\begin{figure}[t!]
	\centering
	\includegraphics[page=1, width=.24\textwidth]{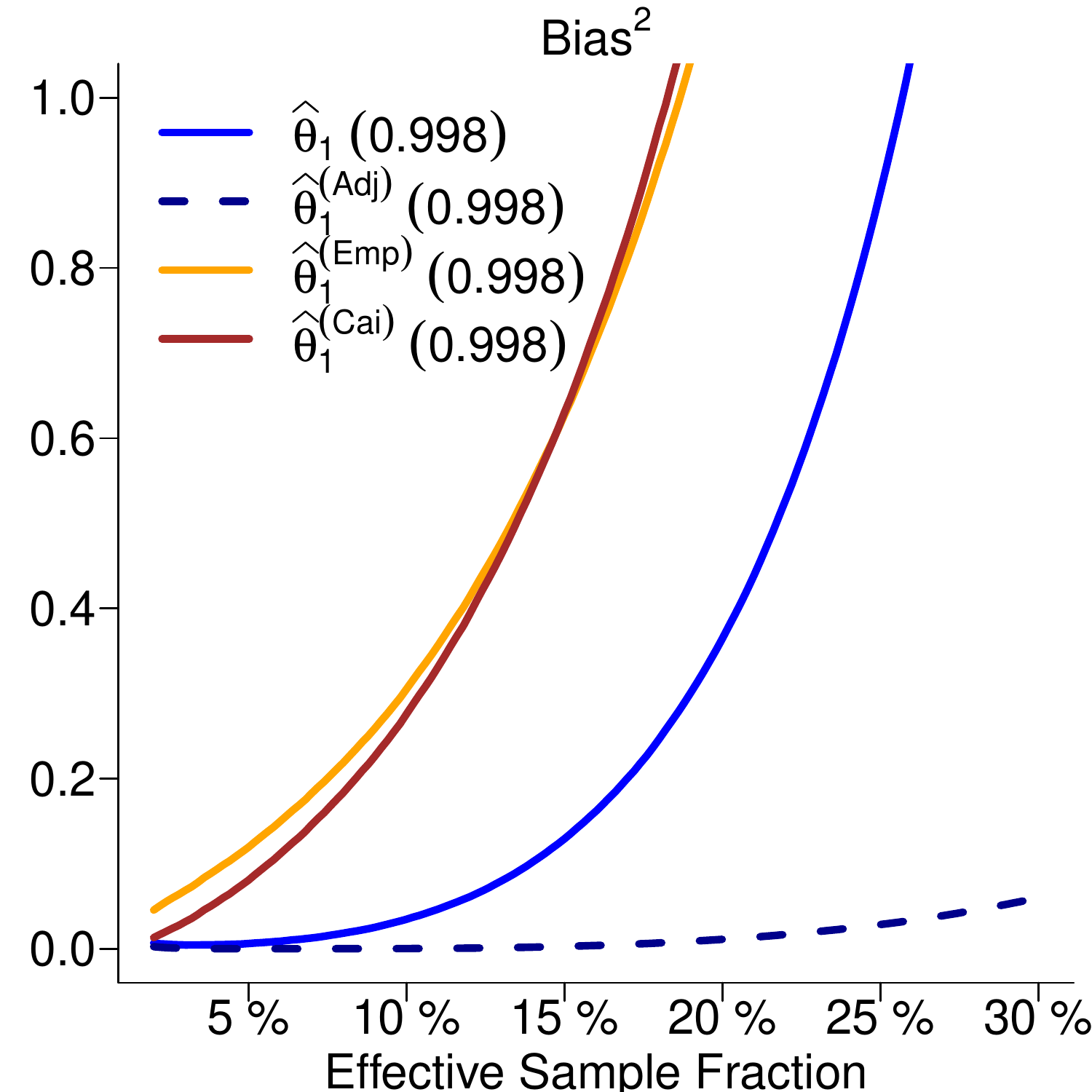}
	\includegraphics[page=2, width=.24\textwidth]{ClaytonHalfT}
	\includegraphics[page=3, width=.24\textwidth]{ClaytonHalfT}
	\includegraphics[page=4, width=.24\textwidth]{ClaytonHalfT}\\
	\includegraphics[page=1, width=.24\textwidth]{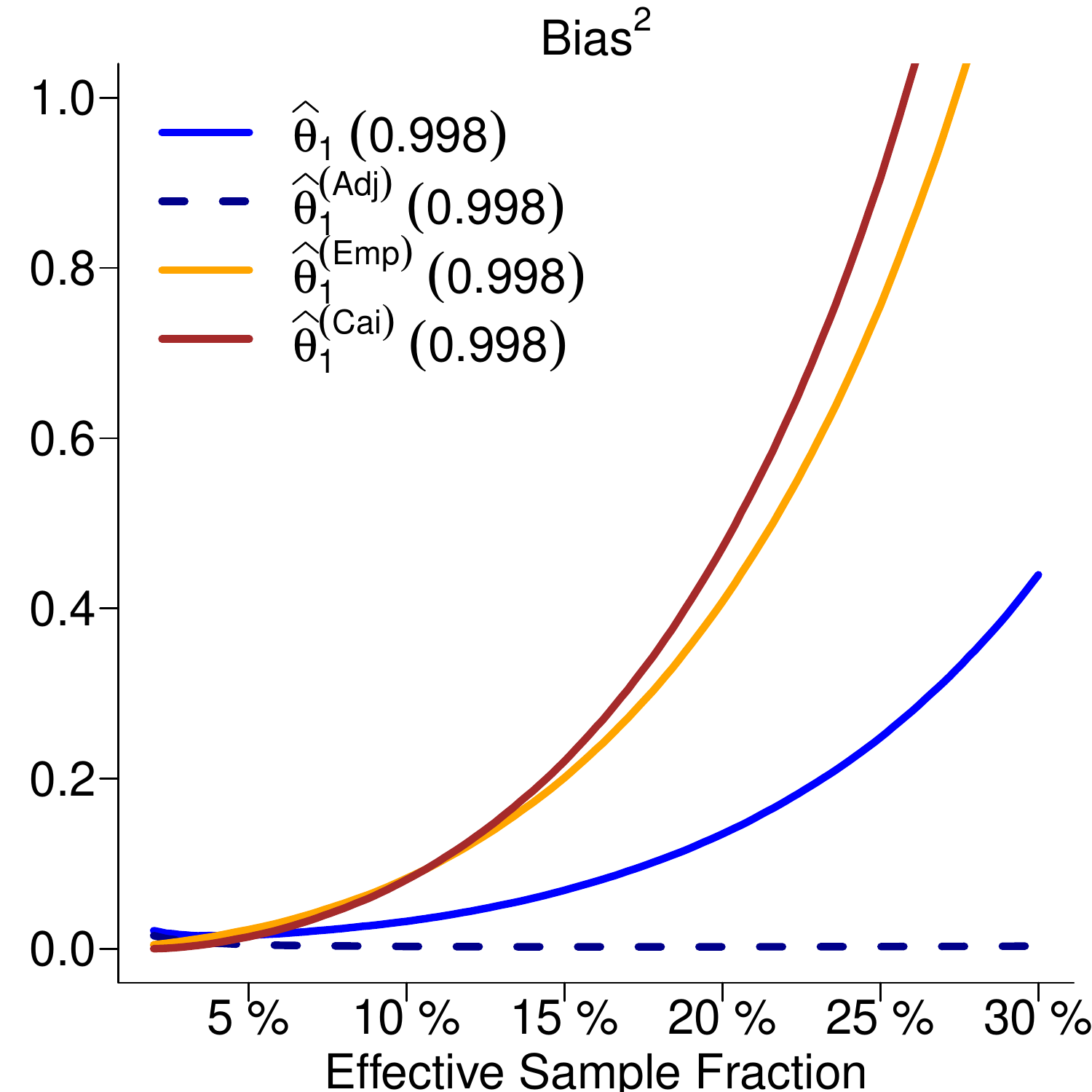}
	\includegraphics[page=2, width=.24\textwidth]{GumbelBurr}
	\includegraphics[page=3, width=.24\textwidth]{GumbelBurr}
	\includegraphics[page=4, width=.24\textwidth]{GumbelBurr}\\
	\includegraphics[page=1, width=.24\textwidth]{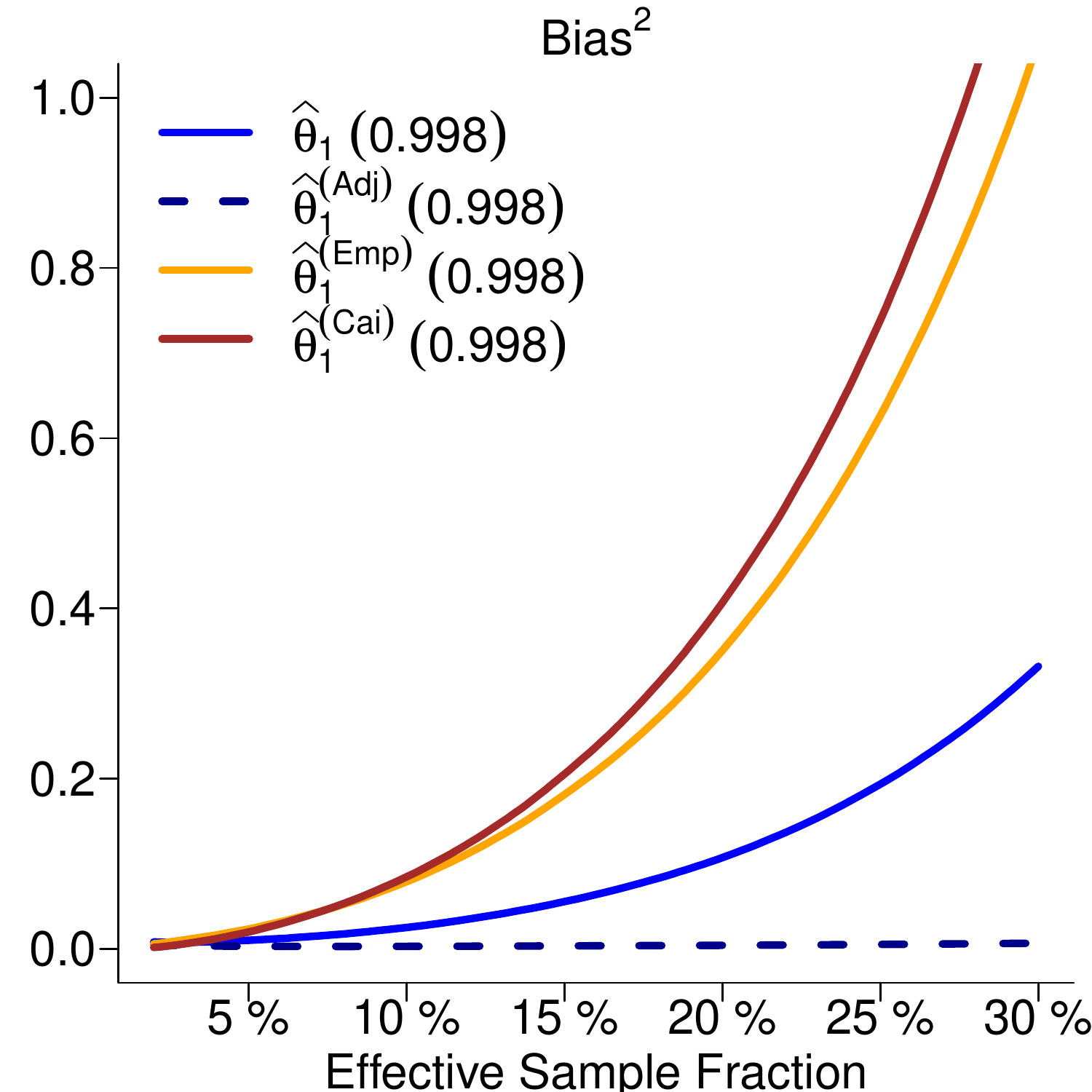}
	\includegraphics[page=2, width=.24\textwidth]{StudentTBurr}
	\includegraphics[page=3, width=.24\textwidth]{StudentTBurr}
	\includegraphics[page=4, width=.24\textwidth]{StudentTBurr}\\
	\includegraphics[page=3, width=.24\textwidth]{4Diff_Burr}
	\includegraphics[page=3, width=.24\textwidth]{4Diff_MAX}
	\includegraphics[page=3, width=.24\textwidth]{4Diff_MEAN}
	\includegraphics[page=4, width=.24\textwidth]{4Diff_Burr}\\
	\includegraphics[page=1, width=.24\textwidth]{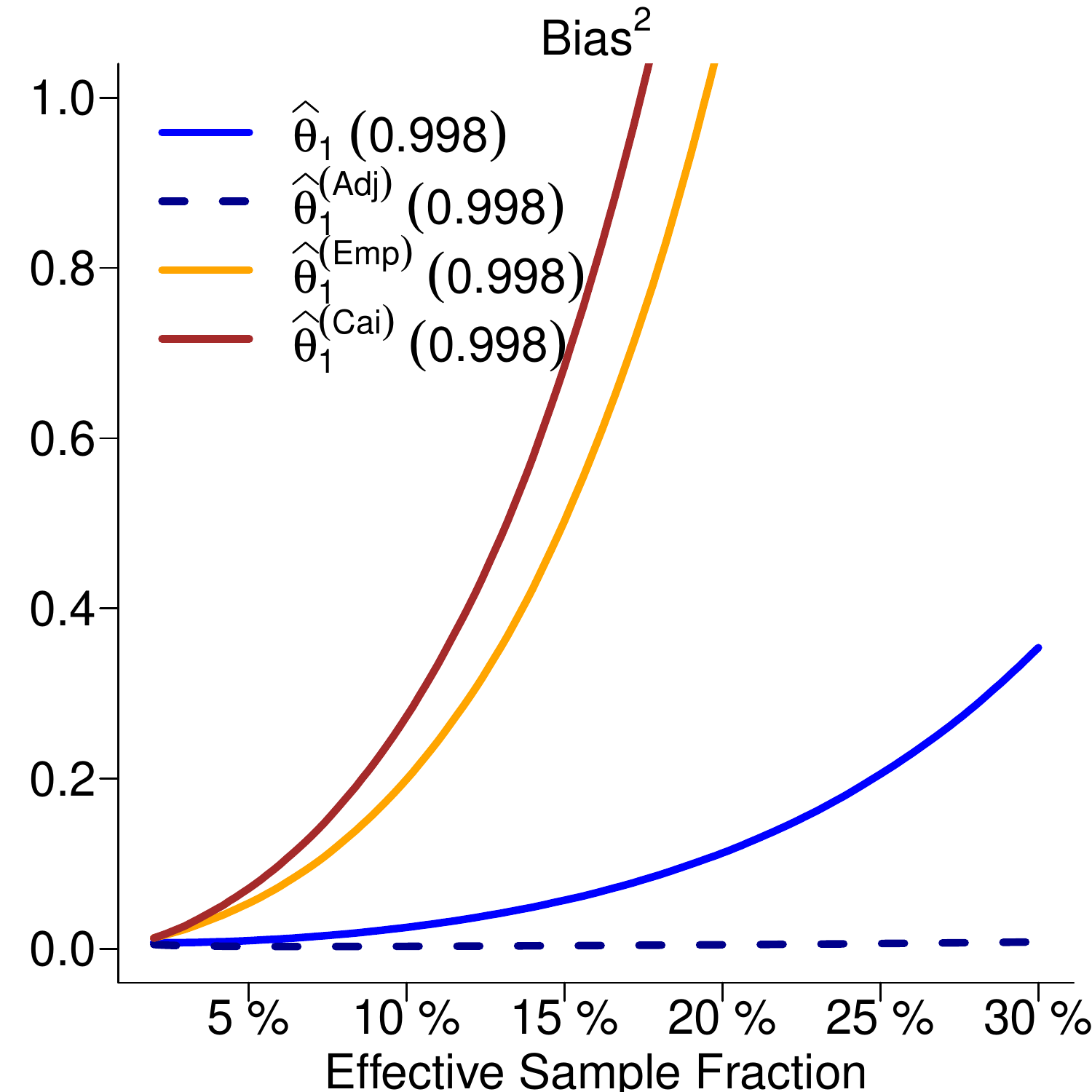}
	\includegraphics[page=2, width=.24\textwidth]{15Dim}
	\includegraphics[page=3, width=.24\textwidth]{15Dim}
	\includegraphics[page=4, width=.24\textwidth]{15Dim}\\
	\caption{Squared-bias, variance, MSE (the first 3 columns from the left) of our estimators (blue solid and dotted lines) and the competitors one (solid brown and yellow solid lines) for different values of $k_n/n$, obtained with the models (i)-(v) (top to bottom row). The very right column reports the non-coverage probability of the four confidence intervals defined in the previous section, with 5\% nominal level.} \label{fig:res_simu_all_dim}
\end{figure}

Figure \ref{fig:res_simu_all_dim} reports the squared-bias, variance, MSE and the coverage probability, from the left to the right column, computed for the effective sample size $k_n$ that ranges between $1\%$ to $30\%$ of $n$. The results obtained with the models (i)-(v) are reported from the top to the bottom row. 
The first three columns from the left,  report the results obtained with our estimator $\widehat{\theta}_{j}$ and its adjusted version $\widehat{\theta}_j^{\text{(Adj)}}$, which are reported with the blue solid and dotted lines, respectively, and the results obtained with the competitors $\widehat{\theta}^{\text{(Emp)}}$ and $\widehat{\theta}_j^{\text{(Cai)}}$, which are reported by the yellow and brown solid lines, respectively.  With model (iv) the marginal distributions are all different, then a partial set of results are reported in the fourth row, to save space. These are, from left to right, the MSE of the MES estimator regarding the Burr marginal distribution and the maximum MSE and the average MSE obtained with respect to the MES estimators of the four marginal distributions. In the very right column the non-coverage probabilities, at 5\% nominal level, concerning the confidence interval defined in \eqref{eq:app_ci_mse} and its adjusted version are reported in yellow and blue dotted lines, respectively, while those concerning the corrected confidence interval defined in \eqref{eq:app_ci_corrected_mse} and its adjusted version are reported in yellow and blue solid lines. The horizontal red dotted line indicates the 5\% nominal level. In the case of model (iv) the results regard the interval for the MES concerning the Burr marginal distribution.

\begin{table}[t!]
\begin{center}
\begin{tabular}{ccccc}
\toprule
Model & $\theta_1(0.998)$ & $\theta_2(0.998)$ & $\theta_3(0.998)$ & $\theta_4(0.998)$\\
\midrule
(i) & $16.58656$ & -- & -- & -- \\
(ii) & $10.09849$ & -- & -- & -- \\
(iii) & $5.270914$ & -- & -- & -- \\
(iv) & $6.965690$ & $3.783465$ & $3.875493$ & $3.869831$ \\
(v) & $6.738795$ & -- & -- & -- \\
\bottomrule
\end{tabular}
\caption{Values of $\theta_j(0.998)$ obtained through intensive Monte Carlo simulations.}
\label{tab:true_val_MES}
\end{center}
\end{table}
According to the squared-bias, variance and MSE, our proposed estimators outperform the competitors with all the considered models. The benefit of our proposals is particularly remarkable in terms of squared-bias and MES and when considering high dimensional models. This result was expected, because unlike the competitors 
our estimator is made up of an estimator of the tail dependence of $F_{\bfX}$ (among those of the tail index and the radial quantile) which significantly improve the resulting MES estimates. The bias corrected version of our estimator shows an excellent performance in terms of MSE as it cuts down considerably the bias term and also the variability.  According to the coverage probability, the interval in formula \eqref{eq:app_ci_mse} performs poorly as the smallest non-coverage probability achieved is approximately $20\%$ with $k_n/n$ between $15\%$ and $20\%$ with all the models. Its adjusted version performs the same for small values of $k_n$, but it does better for large values of $k_n$ (it reaches approximately the nominal level with model (i) when $k_n/n=20\%$  and with models (ii), (iii) and (v) when $k_n/n=30\%$ approximately, while it performs poorly with model (iv), where the smallest non-coverage is obtained approximately with $k_n/n=20\%$).  Lastly, our proposed refined interval in formula \eqref{eq:app_ci_corrected_mse} provides the desired results.  The non-coverage probability achieves the nominal level when $k_n/n\geq10\%$ and remaining stable at such a value, with models (i)-(iii) and (v). In the most difficult case given by model (iv) (all the marginal distributions are different) non-coverage probability achieves approximately the nominal level  when $k_n/n\in[3\%, 18\%]$, while the results slightly deteriorate afterward. The adjusted version of this interval is slightly more conservative than the unadjusted one for $k_n/n\geq15\%$ with models (i), (iii) and (v) and therefore it does not provide a valid alternative.
Concluding, our estimators outperform by far the existing competitors and our best confidence interval is very accurate.

%
\section{Marginal expected shortfall inference for the  macroprudential policy application}\label{sec:realdata}

Macroprudential regulation aims to maintain the stability of the overall financial system. This is done by preventing systemic risks and  implementing then actions to protect the economy from financial crises. The MES is mainly used to rank financial institutions according to their contribution to the market risk. Accordingly, the banks with the greatest impact on the financial system are classified as systemically important financial institutions for which a stricter monitoring and regulation is required \citep[e.g.][]{acharya2013}. The MES also helps to directly quantify the amount of capital that each institution is expected to lose during a downturn and determine thus the amount of additional resources needed to avoid being undercapitalised.

We focus on the weekly negative returns (returns for brevity) from May 5, 2010 to March 31, 2022 concerning 23 financial institutions classified as Global Systemically Important Banks and Domestic Systemically Important Banks from US and Canada, see the left column of Table \ref{tab:MES_results}. The data come from the Center for Research in Security Prices database and are available online at {\tt https://www.crsp.org}.
According to the notation of Section \ref{sec:intro}, the return of a week for the $j$th institution is defined as $X_j=\varphi_j(1-P_t^{(j)}/P_{t-1}^{(j)})$, where $P_{t-1}^{(j)}$ and $P_t^{(j)}$ are the weekly opening and closing prices of the $j$th institution's stock and $\varphi_j=(\text{$j$th institution' s capitalisation} $ $/ \text{market capitalisation})$ at March 31, 2022 is the weight of the $j$th institution compared to the market. Overall there are $n=621$ returns. The autocorrelation plots of the weekly observations indicate that there is almost no autocorrelation. 
The pairwise scatter plots indicate that there is an important dependence between institutions' returns. In this regard, we have computed the pairwise correlations and coefficients of upper tail dependence, whose range is $[0.47,0.87]$ and $[0.25, 0.80]$, respectively, indicating that there is a considerable dependence in the body of the bivariate distributions
but also in their joint upper  tail. We recall that the coefficients of upper tail dependence $\lambda_U$ \citep[e.g.][p. 62]{joe2014}, 
satisfies $0\leq\lambda_U\leq 1$ with the lower and upper bounds representing the case of complete dependence and independence. 
All these results are available in Section D of the supplementary material to save space.
\begin{figure}[t!]
    \centering
    \includegraphics[width = \textwidth/(25)*10]{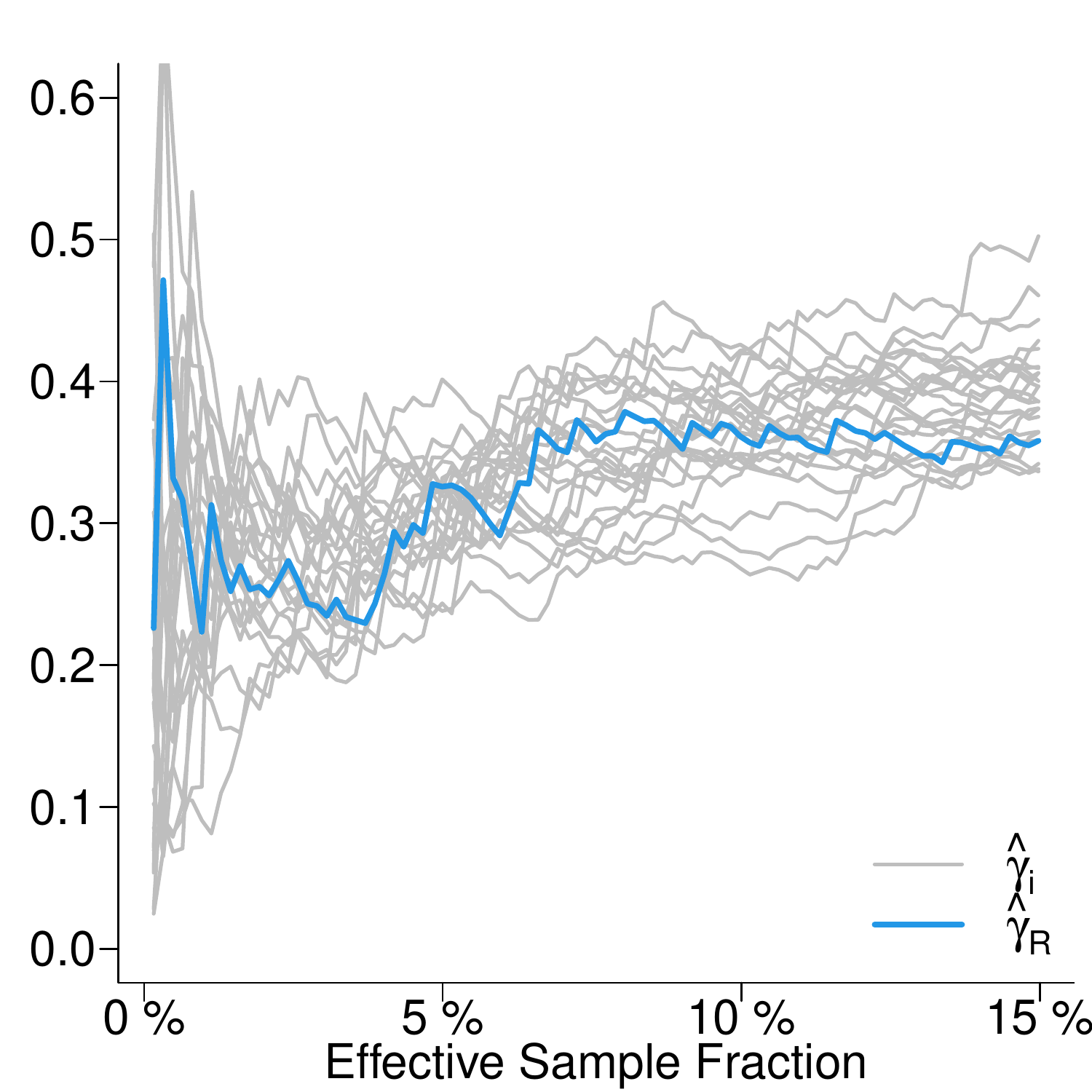}
    \includegraphics[width = \textwidth/(25)*10]{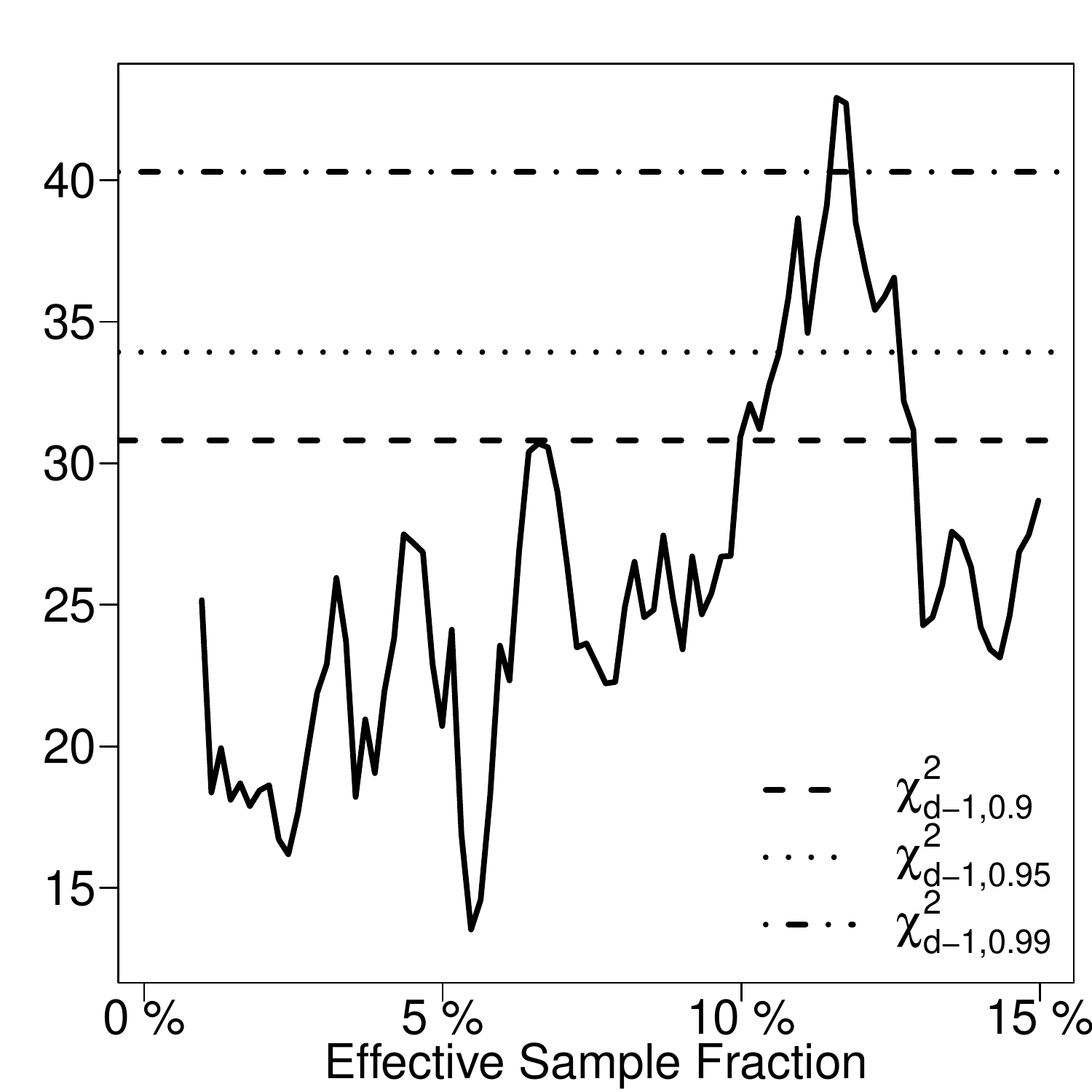}
    \caption{Estimated tail index (left panel) for the individual banks distribution (grey lines) and the weighted sum of them (light-blue line). Test statistic value (solid line in the right panel) to test whether the tail indices of individual banks distributions are equal. Dashed, dotted, dotdash lines are the critical values at 10\%,  5\% and 1\%. Both are computed for the effective sample fraction in $[1\%,15\%]$. }
    \label{fig:tail_index_test}
\end{figure}
The right plot of Figure \ref{fig:tail_index_test} reports the estimated tail index (grey solid lines) for each individual bank and that for the sum of the all them (the light blue solid line), representative of the northern American financial sector. The tail index has been estimated for an effective sample fractions ranging between 1\% to 15\% and its value seems to be stable from the level of about $7\%$ onward. In this region, individual estimates are between about 0.28 to 0.4 and that of the sum is around 0.35. We have performed the statistical test proposed by \citet{padoan2022} to check whether the assumption of equal marginal tail indices is plausible. In particular, the null hypothesis $H_0: \gamma_1=\cdots=\gamma_{23}$ is rejected whenever $\Lambda_n>\chi^2_{22,1-\alpha}$ with significance level $\alpha\in(0,1)$, where $\Lambda_n$ is the test statistic and $\chi^2_{22,1-\alpha}$ is $(1-\alpha)$-quantile of a Chi-square distribution with $22$ degrees of freedom \citep[see][for details]{padoan2022}. The right plot of Figure \ref{fig:tail_index_test} shows the value of the test statistic against the effective sample fraction, which is always below to the critical value at 10\% (dashed line), 5\% (dotted) and 1\% (dotdash) significance level, apart from 9, 6 and 2 single points, respectively. 
All these results support the use of our method for analysing the extreme MES of the financial institutions.
%
\begin{table}[t!]
\centering
\centerline{
\begin{tabular}{lccc}
  \hline
Bank & Size/Rank & MES/Rank & ES/Rank \\
  \hline
JPMorgan Chase \& Co. & 17.60/ 1& 18.88/ 1& 105.10/ 2 \\
  Bank of America Corp. & 14.54/ 2& 17.33/ 2 & 111.78/ 1\\
  Wells Fargo \& Company & ~8.08/ 3& ~8.07/ 3& ~65.12/ 3 \\
  Morgan Stanley & ~6.81/ 5& ~7.75/ 4& ~42.56/ 4\\
  Citigroup Inc. & ~4.62/ 8& ~5.83/ 5& ~39.36/ 5\\
  The Goldman Sachs Group, Inc. & ~4.93/ 7& ~5.10/ 6 & 22.99/ 11\\
  Royal Bank of Canada & ~6.83/ 4 & ~4.60/ 7& ~25.54/ 8\\
  The Toronto-Dominion Bank & ~6.30/ 6& ~4.49/ 8 & ~31.68/ 6 \\
  Truist Financ. Corp. & 3.29/ 13& ~3.48/ 9& ~26.48/ 7\\
  The PNC Financ. Services Group, Inc. & 3.37/ 11& 3.25/ 10& 21.96/ 12\\
  U.S. Bancorp & 3.45/ 10 & 3.16/ 11& 23.80/ 10\\
  The Bank of Nova Scotia & 3.78~/ 9& 2.58/ 12& ~23.23/ 9\\
  Capital One Financ. Corp. & 2.33/ 15 & 2.55/ 13& 18.48/ 13\\
  Bank of Montreal & 3.34/ 12 & 2.54/ 14& 15.52/ 14\\
  The Bank of New York Mellon Corp. & 1.75/ 16 & 1.72/ 15& 12.07/ 16\\
  Canadian Imperial Bank of Comm. & 2.40/ 14 & 1.65/ 16& 15.16/ 15\\
  State Street Corp. & 1.39/ 17 & 1.47/ 17& ~8.16/ 20\\
  Fifth Third Bancorp & 1.29/ 18 & 1.39/ 18& 10.86/ 17\\
  Huntington Bancshares Incorporated & 0.92/ 21& 1.04/ 19& ~9.77/ 18\\
  Northern Trust Corp. & 1.06/ 19 & 0.99/ 20& ~6.74/ 21\\
  M\&T Bank Corp. & 0.96/ 20 & 0.93/ 21& ~9.59/ 19\\
  Comerica Incorporated & 0.52/ 22 & 0.66/ 22& ~4.61/ 22\\
  Zions Bancorporation & 0.43/ 23& 0.53/ 23& ~3.63/ 23\\
   \hline
\end{tabular}}
\caption{List of the main northern American banks and their: market capitalisation (size in \%), marginal expected shortfall (in \%) and expected shortfall (in billions of \$).}
\label{tab:MES_results}
\end{table}
Table \ref{tab:MES_results} reports the size of a bank (i.e. $\varphi_j\cdot 100\%$), and the MES computed with our approach (both in percentage) and the ES in billions of USD for comparison (from left to right). On their basis the ranking of banks from the largest to the smallest value are also reported. According to the MES (the contribution of the $j$th institution to the market downturn), four institutions, namely JP Morgan, Bank of America, Wells Fargo and Morgan Stanley are responsible for more than 50\% of the northern American market downturn. This ranking is not quite the same as that obtained with the size of the bank, suggesting that the market capitalisation alone does not fully explain the overall capital shortfall, which is  influenced by other factors. Finally, Table \ref{tab:MES_results} reports the estimate of the expected capital loss of the $j$th institution (in percentage), i.e. $\varphi_j^{-1}\theta_j(\tau)\cdot 100\%$ and the corresponding $95\%$ confidence interval, with $\tau=0.9989$ and $\tau=0.9995$ that are such that $R$ is expected to exceed $Q_R(\tau)$, on average every $10$ and $20$ years. Interestingly to see that those who contribute the most to  the northern American market downturn are not necessarily those who lose the majority of their individual capital.
\begin{table}[ht]
\centering
\centerline{
\begin{tabular}{lcccc}
\hline
 Bank   & \multicolumn{2}{c}{10 years} & \multicolumn{2}{c}{20 years}  \\
        & MES & 95\%CI & MES & 95\%CI \\
  \hline
Comerica Incorporated & 32.53 & [19.56, 54.10] & 41.84 & [23.55, 74.37] \\
  Citigroup Inc. & 32.05 & [19.27, 53.30] & 41.22 & [23.19, 73.26] \\
  Zions Bancorporation & 30.86 & [18.56, 51.33] & 39.70 & [22.34, 70.55] \\
  Bank of America Corp. & 30.32 & [18.23, 50.43] & 39.01 & [21.95, 69.32] \\
  Morgan Stanley & 28.97 & [17.42, 48.18] & 37.26 & [20.97, 66.23]\\
  Huntington Bancshares Incorporated & 28.82 & [17.33, 47.93] & 37.07 & [20.86,  65.89] \\
  Capital One Financ. Corp. & 27.91 & [16.78, 46.42] & 35.90 & [20.20, 63.80]\\
  Fifth Third Bancorp & 27.32 & [16.42, 45.43] & 35.14 & [19.77, 62.45]\\
  JPMorgan Chase \& Co. & 27.29 & [16.41, 45.39] & 35.11 & [19.75, 62.40]\\
  Truist Financ. Corp. & 26.91 & [16.18, 44.76] & 34.62 & [19.48, 61.52]\\
  State Street Corp. & 26.83 & [16.13, 44.62] & 34.51 & [19.42, 61.33] \\
  The Goldman Sachs Group, Inc. & 26.32 & [15.83, 43.78] & 33.86 & [19.05, 60.18]\\
  Wells Fargo \& Company & 25.39 & [15.27, 42.22] & 32.66 & [18.38, 58.04]\\
  The Bank of New York Mellon Corp. & 24.94 & [15.00, 41.48] & 32.08 & [18.05,  57.01]\\
  M\&T Bank Corp. & 24.78 & [14.90, 41.22] & 31.88 & [17.94, 56.66] \\
  The PNC Financ. Services Group, Inc. & 24.53 & [14.75, 40.80] & 31.56 & [17.76, 56.09]\\
  Northern Trust Corp. & 23.71 & [14.26, 39.44] & 30.50 & [17.16, 54.21]\\
  U.S. Bancorp & 23.30 & [14.01, 38.75] & 29.97 & [16.86, 53.26]\\
  Bank of Montreal & 19.35 & [11.64, 32.18] & 24.89 & [14.01, 44.24]\\
  The Toronto-Dominion Bank & 18.14 & [10.90, 30.16] & 23.33 & [13.13, 41.46]\\
  Canadian Imperial Bank of Commerce & 17.53 & [10.54, 29.15] & 22.55 & [12.69,  40.07]\\
  The Bank of Nova Scotia & 17.37 & [10.44, 28.88] & 22.34 & [12.57, 39.70]\\
  Royal Bank of Canada & 17.13 & [10.30, 28.48] & 22.03 & [12.40, 39.15]\\
   \hline
\end{tabular}}
\caption{List of the main northern American banks and their expected capital loss.}
\label{tab:MES_results}
\end{table}

%
%
\section{Extension to serially dependent data}\label{sec:extension}

In the financial industry certain applications require to deal with high-frequency data, and to avoid erroneous and overly-optimistic estimation results, their serial dependence cannot be disregarded. For this purpose we extend the theory developed for i.i.d. data to the case of time dependent data. To cover important and popular financial time series models such as ARMA and GARCH (among other) we focus on a multivariate $\beta$-mixing heavy-tailed strictly stationary time series framework, see \citet{hoga2018}, and see for example \cite{drees2003} and \citet{davison2022tail} for related results in the univariate and bivariate framework.  

We consider the following setup. Let $(\bfX_i, i=1,2,\ldots)$ be a multivariate strictly stationary time series with a continuous joint marginal distribution $F$, where $\bfX_i=(X_{i,1},\ldots,X_{i,d})^\top$. 
\begin{cond}\label{cond:multi_ts_cond}
The time series $(\bfX_i, i=1,2,\ldots)$ satisfies the following properties:
\begin{enumerate}
\item[(i)] it is $\beta$-mixing, i.e. for any $m\geq1$, let $\mathcal{F}_{1,m}=\sigma(\bfX_1,\ldots,\bfX_m)$ and $\mathcal{F}_{m,\infty}=\sigma(\bfX_m,\bfX_{m+1},\ldots)$ be the past and future sigma-fields generated by the time series, then 
$$
\beta(l):=\sup_{m\geq1}\Expect\left(\sup\{|\Prob(B|\mathcal{F}_{1,m})-\Prob(B)|:B\in\mathcal{F}_{m+l,\infty}\}\right)\to0,\quad l\to\infty;
$$
\item[(ii)] its marginal distribution $F$ satisfies 
Condition \ref{cond:mvt_sec_order};
\item[(iii)] there are integer sequences $(l_n)$ and $(r_n)$ such that $l_n\to\infty$, $r_n\to\infty$, $n\beta(l_n)/r_n\to 0$ and $r_nk_n^{-1/2}\log^2k_n\to 0$ as $n\to \infty$, with $l_n=o(r_n)$ and $r_n=o(n)$.
\item[(iv)] for any $t=1,2\ldots$ there are functions $r_t$ on $[0,\infty]^2\setminus\{(\infty,\infty)\}$ such that
$$
\lim_{s\to\infty}s\Prob\left(\overline{F}_{R_1}(R_1)\leq\frac{x}{s},\overline{F}_{R_{1+t}}(R_{1+t})\leq\frac{y}{s}\right)=r_t(x,y),\;\forall (x,y)\in[0,\infty]^2\setminus\{(\infty,\infty)\};
$$
\item[(v)] there are nonnegative $\eta(t)$ such that $\sum_{t=1}^\infty\eta(t)<\infty$ and $D\geq 0$ satisfying for all $t=1,2\ldots$ and $x_1,x_2,y_1,y_2\in[0,1]$ with $x_1<x_2$ and $y_1<y_2$ as $s\to\infty$
\begin{align*}
s\Prob\left(\frac{x_1}{s}<\overline{F}_{R_1}(R_1)\leq \frac{x_1}{s},\frac{y_1}{s}<\overline{F}_{R_{1+t}}(R_{1+t})\leq \frac{y_1}{s}\right)&\leq\eta(t)\sqrt{(x_1-x_2)(y_1-y_2)}\\
&+\frac{D}{s}(x_1-x_2)(y_1-y_2).
\end{align*}
\end{enumerate}
Condition \ref{cond:multi_ts_cond}(i) is the multivariate $\beta$-mixing condition adopted by \citet{hoga2018} that extends in a natural way the univariate one \cite[see e.g.][]{drees2003}. Condition \ref{cond:multi_ts_cond}(ii) is the same used in the independence case. Conditions \ref{cond:multi_ts_cond}(iii)-\ref{cond:multi_ts_cond}(v) are standard conditions adopted in extreme value analysis with mixing conditions \citep[see][]{davison2022tail,drees2003}, which are required for the radial component, in our context. The sequences $(r_n)$ and $(l_n)$ are those used in the big-blocks separated by small-blocks approach in the literature on mixing time series, where the sequence $(\bfX_i, i=1,2\ldots)$ is split into $m_n=\round{n/(r_n+l_n)}$  big-blocks of size $r_n$ separated by small-blocks of size $l_n$ (and a negligible remainder). We now provide a general result for serially dependent data.
\end{cond}
\begin{theo}\label{theo:asym_multi_ts}
Work under Condition \ref{cond:multi_ts_cond}. Assume also that 
$\sqrt{k_n}A_R(n/k_n)\to\lambda_R\in \Real$, as $n\to\infty$. Then, for any $j\in\{1,\ldots,d\}$, as $n\to\infty$
$$
\frac{\sqrt{k_n}}{\log(k_n/n(1-\tau_n))}\log \frac{\widehat{\theta}_j(\tau_n)}{\theta_j(\tau_n)}
	\stackrel{d}{\to}
	\mathcal{N}\left(\frac{\lambda_R}{1- \rho_R},\gamma^2\left(1+2\sum_{t=1}^\infty r_{t}(1,1)\right)\right).
$$
\end{theo}
As in \citet{davison2022tail,drees2003}, the sequence of functions $(r_t, t=1,2\ldots)$ capture the extremal dependence within the time series between time instants. In the i.i.d. case we have that $r_t=0$ for all $t\geq1$ and the above asymptotic variance simplifies to that in Theorem \ref{theo:est}.
\section*{Acknowledgements}
Simone Padoan is supported by the Bocconi Institute for Data Science and Analytics (BIDSA), Italy. 

\bibliographystyle{chicago} 
\bibliography{bibliopm_final}

\end{document}